\documentclass{article}

\usepackage{amssymb,amscd,amsthm,latexsym}
\usepackage[all]{xy}

\newtheorem{defi}{Definition}[section]
\newtheorem{theo}{Theorem}[section]
\newtheorem{lemma}{Lemma}[section]
\newtheorem{prop}{Proposition}[section]
\newtheorem{coro}{Corollary}[section]

\def\into{ \rightarrowtail }

\def\splito{ \rightleftarrows }

\newdir{ >}{{}*!/-8pt/@{>}}
\def\EE{ \mathbb{E} }
\def\CC{ \mathbb{C} }
\def\DD{ \mathbb{D} }
\def\VV{ \mathbb{V} }

\def\Pt{ \mathrm{Pt} }

\newcommand{\Ker}{ \ensuremath{\mathrm{Ker}} }

\newcommand{\Eq}{ \ensuremath{\mathrm{Equ}}}

\newcommand{\SkB}{\mathsf{SkB}}
\newcommand{\Cat}{\mathsf{Cat}}
\newcommand{\Grd}{\mathsf{Grd}}
\newcommand{\RGh}{\mathsf{RGh}}
\newcommand{\SRg}{\mathsf{SRg}}
\newcommand{\Rng}{\mathsf{Rng}}

\newcommand{\Gp}{\mathsf{Gp}}

\newcommand{\Set}{\mathsf{Set}}

\newcommand{\Mon}{\mathsf{Mon}}

\newcommand{\Aff}{\mathsf{Aff}}

\newcommand{\Sub}{\mathsf{Sub}}

\newcommand{\Norm}{\mathsf{Norm}}

\begin{document}

\author{Dominique Bourn}

\title{$\Eq$-saturating categories}

\date{}

\maketitle

\begin{abstract}
Starting from the varietal notion of syntactic equivalence relation, we generalized it to a categorical concept; namely $\Eq$-saturating category. We produce various examples and focuse our attention on the protomodular context in which any equivalence relation is then shown to have a centralizer.	
\end{abstract}

\section*{Introduction}

In a finitely complete category $\EE$ a monomorphism $u: U\into X$ is normal to an internal equivalence relation $R$ on $X$ when we get $u^{-1}R=\nabla_U$ and when moreover, in the induced diagram:
$$
\xymatrix@=20pt
{  U\times U \ar@<-1ex>[d]_{p_0^U}\ar@<+1ex>[d]^{p_1^U}\ar@{ >->}[r]^{\bar{u}}
	& R\ar@<-1ex>[d]_{d_0^R}\ar@<+1ex>[d]^{d_1^R}
	\\
	U\ar[u]\ar@{ >->}[r]_u & X \ar[u]
}
$$
the square indexed by $0$ (or equivalently by $1$) is a pullback.
In any protomodular category \cite{Bprot}, a monomorphism $u$ is normal to at most one equivalence relation \cite{Bnorm}, as it is the case in the category $\Gp$ of groups. See also \cite{BM} for many other aspects of normal subobjects in a category $\EE$.

\medskip

In the category $\Mon$ of monoids, the notion of {\em syntactic equivalence relation} allows to show that, when it is non empty, the preordered set of equivalence relations $R$ w.r. to which a monomorphism $u$ is normal has a supremum \cite{schu} (back to 1956!).

\medskip

The aim of this work is to investigate the categories in which such a kind of {\em generalized syntactic equivalence relation} and such a kind of supremum do exist. We call them $\Eq\EE$-saturating categories, see precise Definition \ref{supnorm}. Following \cite{Alm}, it is the case of any variety $\VV$ of Universal Algebra. We show that it is the case as well for any fiber $\Cat_X\EE$ of the fibration $(\;)_0: \Cat\EE \to \EE$ (where $\Cat\EE$ is the category of internal categories in $\EE$) when $\EE$ is $\Set$ or when $\EE$ is any Mal'tsev category in the sense of \cite{CLP}, \cite{CPP}  or any $\Eq\EE$-saturating Gumm category in the sense of \cite{BG}; the three results following from quite heterogeneous arguments.

\medskip

When a category $\EE$ is both protomodular and  $\Eq\EE$-saturating, we show that any equivalence relation $R$ on an object $X$ has necessarily a centralizer. The protomodular varieties were characterized in \cite{BJ}. The exact protomodular context allows to compare with precision the $\Eq$-saturating axiom with the existence of abstract normalizers in the sense of \cite{Gr}, see Section \ref{normali}. This comparison is of a complementary nature. An additive category producing a particular case of protomodularity, we incidentally show that any additive category is trivially $\Eq$-saturating.

\medskip

This article is organized along the following lines. Section 1 is devoted to notations and Section 2  to recalls about the syntactic equivalence relation. In Section 3, the notion of $\Eq$-saturating category is introduced. Section 4 is devoted to $\Eq$-saturating protomodular categories and Section 5 to the fibers $\Cat_X\EE$. I thank very much J-\'E. Pin for his help in my search of precise references.

\section{Notations}

Any category $\EE$ will be supposed finitely complete. Recall that an internal relation $S$ on an object $X$ is a subobject of the product $X\times X$:
$$
\xymatrix@=30pt 
{S\; \ar@{>->}[r]^-{(d_0^S,d_1^S)}
	&X\times X}
$$
We shall use the simplicial notation for describing any internal reflexive or equivalence relations on $X$:
$$
\xymatrix@=30pt
{S\ar@<+2ex>[r]^-{d^S_1}\ar@<-2ex>[r]_-{d^S_0}
	&X\ar[l]|{s^S_0}
}
$$

We denote by $\Eq\EE$ the category of internal equivalence relations in $\EE$ whose morphisms are given by the pairs $(f,\bar f)$ making the following diagram commute:
$$
\xymatrix@=25pt
{ S\ar@<-1ex>[d]_{d_0^S}\ar@<+1ex>[d]^{d_1^S}\ar[r]^{\bar{f}}
	&R\ar@<-1ex>[d]_{d_0^R}\ar@<+1ex>[d]^{d_1^R}
	\\
	X\ar[u]\ar[r]_f &Y\ar[u]
}
$$
Since $\bar f$, when it exists, is uniquely determined by $f$, we shall denote  it $f:S\to R$.
We denote by $(\;)_0: \Eq\EE\to \EE$ the forgetful functor associating with any equivalence relation its ground object; as soon as $\EE$ has finite  limits, this functor $(\;)_0$ is a fibration whose {\em cartesian maps} are given by the inverse images $f: f^{-1}R\to R$ along the ground morphism $f: X\to Y$ in $\EE$. Each fiber $\Eq_X\EE$ of this fibration is a preorder, which is equivalent to the fact that the functor $(\;)_0$ is faithful. 

We classically denote by $\Delta_X$ the discrete equivalence on the object $X$ which is the smallest object in the fiber $\Eq_X\EE$ above $X$ and by $\nabla_X$ the undiscrete equivalence relation which is the largest one. The \emph{kernel equivalence relation} $R[f]$ of a map $f:X\to Y$ is then the domain of the cartesian map above $f$ with codomain $\Delta_Y$, it is given by the following pullback in $\Eq\EE$:
\[\xymatrix@C=2pc@R=1,5pc{ R[f]  \ar@{ >->}[d]_{1_X} \ar[r]^{f} & \Delta_Y \ar@{ >->}[d]^{1_Y}\\
	\nabla_X \ar[r]_{\nabla_f} & \nabla_Y}\]

A morphism $f: S\to R$ of equivalence relations is {\em fibrant} when the square indexed by $0$ (or equivalently  the square indexed by $1$) is a pullback.
Fibrant maps and cartesian maps are stable under composition, under pullback along any morphism in $\Eq\EE$, and such that when $g.f$ and $f$ belong to any of these classes, so is $f$. Any fibrant morphism above a monomorphism u: $U\into X$ is necessarily cartesian, namely an inverse image. Such a fibrant monomorphism $u: S\into R$ of equivalence relations is said to be a {\em saturated} one. Any monomorphism $\Delta_U: \Delta_X\into \Delta_Y$ is saturated. When  $u$ and $v$ are monomorphisms and when the composition $S\stackrel{u}{\into} R \stackrel{v}{\into} T$ is saturated, then so is $S \stackrel{u}{\into}R$. 
\begin{defi}
	A monomorphism $u: U\into X$ is said to be saturated  w.r. to an equivalence relation $R$ on $X$ when the cartesian map above $u$ with codomain $R$:
	$$ \xymatrix@=25pt{
		 {u^{-1}R\;} \ar@{>->}[r]^{\tilde{u}} \ar@<-2ex>[d]_{p_0^U} \ar@<2ex>[d]^{p_1^U} & R   \ar@<-2ex>[d]_{d_0^R}\ar@<2ex>[d]^{d_1^R}\\
		{U\;} \ar@{>->}[r]_{u} \ar[u]|{s_0^U} & X \ar[u]|{s_0^R}
	}   
	$$
	is saturated.
\end{defi}
So, in any category $\EE$, a subobject $u: U\into X$  is {\em normal to} $R$ if and only if it is saturated w.r. to $R$ and such that $u^{-1}R=\nabla_U$.
In $\Set$, when $U$ is a non-empty part of $X$, the inclusion $u: U\into X$ is saturated w.r. to $R$ when:
$ \forall (x,x')\in X\times X, \; xRx' \Rightarrow [x\in U \iff x'\in U] $. 
In other words, this inclusion is saturated w.r. to $R$ if and only if it is a union of equivalence classes of $R$.
\begin{lemma}\label{lemma0}
Given any category $\EE$ and any monormorphism $u: U\into X$,\\
1) if $u$ is saturated w.r. to $R$ and $S\subset R$, then  $u$ is saturated w.r. to $S$;\\
2) if $u$ is normal to $R$, then $u$ is saturated w.r. to any equivalence relation $S\subset R$;\\
3) when the following square is a pullback (i.e. when $u=f^{-1}(v)$):\\
\[\xymatrix@C=2pc@R=1pc{ {U\;} \ar[d]_{f'} \ar@{ >->}[r]^u  & X \ar[d]^{f}\\
	{V\;} \ar@{>->}[r]_{v} & Y}\]
then $u: R[f']\into R[f]$ is saturated and so, $u$ is saturated w.r. to $R[f]$.
\end{lemma}
 
\section{Syntactic equivalence relations}

\subsection{Monoids}

We shall denote  by $\Mon$ the category of monoids.  
Let $M$ be a monoid and $W$  any non-empty part of $M$. Recall from \cite{schu} the following definition and results:
\begin{defi}
	Call {\em syntactic relation} associated with $W$ the relation $R_W$ defined by: $$mR_Wn \iff \forall (x,y)\in M\times M,\;  [xmy\in W \iff xny\in W]$$
\end{defi}
\begin{prop}\label{prop1}
	Given  any non-empty $W$, the relation $R_W$ is an internal equivalence relation in $\Mon$ whose class of the unit $1$ is the monoid $\overline W$ defined by:
	$$\overline W=\{ m\in M/ \forall (x,y)\in M\times M,\; [xmy\in W \iff xy\in W] \}$$
\end{prop}
\proof
The relation $R_W$ is an equivalence relation since it is the case of the hyperrelation [$\iff$]. Let us check it is internal in $\Mon$:\\
1) We get $1R_W1$ since $R_W$ is reflexive;\\
2) Suppose $mR_Wn$ and $m'R_Wn'$. Then we get:\\
$\forall (x,y)\in M\times M,\; xmm'y\in W \iff xnm'y$ since $mR_Wn$, and:\\
$\forall (x,y)\in M\times M,\; xnm'y\in W \iff xnn'y$ since $m'R_Wn'$.\\
In this way, we get $mnR_Wm'n'$.
\endproof
\begin{lemma}
	When we have $1\in W$, we then get: $\overline W\subset W$.
\end{lemma}
\proof
Suppose $m\in \overline W$ taking $(x,y)=(1,1)$, we get: $m\in W \iff 1\in W$.
\endproof
So, when $K$ is a submonoid of $M$, we get $\overline K\subset K$  in $\Mon$. We get $K=\overline K$ under the following condition: 
\begin{defi}
	Let $K$ be a submonoid of $M$. Call it normal when:\\ $\forall k\in K \; {\rm and}\; \forall (x,y)\in M\times M, \; [xky\in K \iff xy\in K]$.
\end{defi}
This definition gives way to a characteristic definition of a kernel inclusion in $\Mon$:
\begin{prop}
	A submonoid $K$ of $M$ is normal in the previous sense if and only if it is a kernel inclusion in $\Mon$.
\end{prop}
\proof
Checking that the kernel of a monoid homomorphism is normal in this sense is straightforward. Conversely, let $K$ be a normal submonoid of $M$. Consider the monoid homomorphism $q: M\to M/R_K$; then $k\in \Ker f$ if and only if $kR_K1$ if and only if $k\in \overline K=K$.
\endproof
\begin{theo}
	Let $K$ be a normal submonoid of $M$. Then $R_K$ is the largest equivalence relation $S$ on $M$ in $\Mon$ to which  $K$ is normal.
\end{theo}
\proof
Let $S$ be an equivalence relation on $M$ in $\Mon$ such that the class of the unit is $K$. So, we have $xmy \in K\iff xmy S1$. Now, $mSn$ implies $xmySxny,\; \forall (x,y) \in M\times M$. So, $\forall (x,y), \; [xmy \in K \iff xny \in K]$, and we get $S\subset R_K$.
\endproof
We have not yet reach the {\em universal property} of the syntactic equivalence relation $R_K$; for that we need the saturated monomorphisms between equivalence relations.
\begin{theo}
	When $K$ is a submonoid of $M$, the internal equivalence relation $R_K$ on $M$ is such that:\\
	1) the inclusion $k: K\into M$ is saturated w.r. to $R_K$;\\
	2) it is the largest internal equivalence relation in $\Mon$ satisfying the saturation property above the inclusion $k$. 
\end{theo}
\proof
1) Suppose $mR_Kn$. Setting $(x,y)=(1,1)$, we get $m\in K \iff n\in K$.\\
2) Suppose $S$ is an equivalence relation w.r. to which $k$ is saturated. From $mSn$, $\forall (x,y)\in X\times X$, we get $xmySxny$. Thus $xmy\in K$ if and only if $xny\in K$; so we get $mR_Kn$. 
\endproof

\subsection{Groups}

\begin{lemma}
	Let $G$ be a group and $K\hookrightarrow G$ a submonoid. Then $K$ is a normal subobject in $\Mon$ if and only if it is a normal subgroup of $G$.
\end{lemma}
\proof
A normal subgroup is a normal submonoid since $\forall k\in K$ and $\forall (x,y)\in G\times G$, we get:
$$ xky=xkx^{-1}xy \in K \iff xy\in K$$
Conversely, suppose $K$ is a normal submonoid. Let $k\in K$. Then with $(1,k^{-1})\in G\times G$, we get: 
$1kk^{-1}=1\in K \iff 1k^{-1}=k^{-1}\in K$. And $K$ is a subgroup of $G$. Now, with any $(x,x^{-1})\in G\times G$, we get $xkx^{-1}\in K \iff xx^{-1}=1\in K$. So, $K$ is a normal subgroup.
\endproof
Let $\Gp$ be the category of groups, $G$ a group, and $K$ a normal subgroup of $G$. We know that the set of internal equivalence relations in $\Gp$ to which $K$ is normal is reduced to one element; accordingly the category $\Gp$ is an extremal example of a category in which the set of internal equivalence relations to which a subobject is normal, when it is non-empty, has a supremum. In $\Gp$, this equivalence relation coincides with the syntactic equivalence relation $R_K$ in $\Mon$. This is a corollary of the following well known
\begin{lemma}\label{lemma2}
	Let $G$ be a group and $S$ an equivalence relation on $G$ which internal in $\Mon$. Then it  is internal in $\Gp$.
\end{lemma}
\proof
It remains to check that when we have $mSn$, we get $m^{-1}Sn^{-1}$. From $mSn$, we get $1=mm^{-1}Snm^{-1}$; then $n^{-1}Sm^{-1}$.
\endproof

\subsection{The varietal case}

The definition of a syntactic equivalence relation with its characteristic universal property can be extended to any variety $\VV$, see \cite{Alm}.
\begin{defi}
	Given any algebra $A$ and any non-empty subset $L$, call {\em syntactic relation} associated with $L$ the relation $R_L$ defined by $mR_Ln$ when:\\ for any $(a_1,\cdots, a_n)\in A^n$ and any term $\tau (x_0,x_1,\cdots, x_n)$ of $\VV$, we get:\\
	$\tau (m,a_1,\cdots, a_n)\in L \iff \tau (n,a_1,\cdots, a_n)\in L$
\end{defi}
\begin{theo}\label{alme}
	When $L$ is a subalgebra of $A$, the relation $R_L$ is a congruence on $A$ in $\VV$ is such that:\\
	1) the inclusion $l: L\into A$ is saturated w.r. to $R_L$;\\
	2) it is the largest congruence in $\VV$ with the saturation property w.r. to $L$. 
\end{theo}

\subsubsection{Semirings}

We shall briefly review the case of semirings, for two reasons:\\
1) it enlarges the strong structural parallelism with monoids initiated in \cite{BMMS},\\
2) the list of axioms of the syntactic equivalence relation can be reduced to a finite one. 

A semiring is a commutative monoid $(A,+,0)$ endowed with an associative binary operation $\cdot: A\times A\to A$ which is distributive w.r. to $+$ and such that $x.0=0$ for all $x\in A$. We shall denote  by $\SRg$ the category of semirings. 
Let $A$ be a semiring and $W$  any non-empty part of $A$. 
\begin{defi}
	Call {\em syntactic relation} associated with $W$ the relation defined by: $mR_Wn$ when the following conditions hold:\\
	1) $\forall x\in A,\; [x+m\in W \iff x+n\in W]$;\\
	2) $\forall (x,y)\in A\times A,\; [x+my\in W \iff x+ny\in W]$;\\
	3) $\forall (x,y)\in A\times A,\; [x+ym\in W \iff x+yn\in W]$;\\
	4) $\forall (x,y,z)\in A\times A\times A,\; [x+ymz\in W \iff x+ynz\in W]$.
\end{defi}
\begin{prop}
	Given  any $W$, the relation $R_W$ is a congruence in $\SRg$ whose class of the unit $0$ is the subsemiring $\overline W$ defined by the set of objects $m\in A$ such that:\\
	1) $\forall x\in A,\; [x+m\in W \iff x\in W]$;\\
	2) $\forall (x,y)\in A\times A,\; [x+my\in W \iff x\in W]$;\\
	3) $\forall (x,y)\in A\times A,\; [x+ym\in W \iff x\in W]$;\\
	4) $\forall (x,y,z)\in A\times A\times A,\; [x+ymz\in W \iff x\in W]$.
\end{prop}
\proof
The relation $R_W$ is an equivalence relation since it is the case of the hyperrelation [$\iff$]. Let us check it is a congruence in $\SRg$.\\ 
Compatibility with the law $+$:\\
1) We get $0R_W0$ since $R_W$ is reflexive.\\
2) Suppose $mR_Wn$ and $m'R_Wn'$. We have to check $(m+m')R_W(n+n')$. Condition 1) holds by Proposition \ref{prop1}.\\
Condition 2): suppose $x+(m+m')y=x+my+my'$ in $W$. Since $mR_Wn$ this the case if and only if $x+ny+m'y\in W$ by Condition 2 on $R_W$. This last point holds if and only if $x+ny+n'y=x+(n+n')y\in W$ since $m'R_Wn'$ by the same condition 2) on $R_W$. \\
Condition 3) is satisfied on the same way by Condition 3) on $R_W$. Same thing for Condition 4).

\noindent Compatibility with the law $\cdot$:\\
Suppose $mR_Wn$ and $m'R_Wn'$. We have to check $mm'R_Wnn'$.\\
Condition 1). We get $x+mm'\in W$ if and only if $x+nm'\in W$ by Condition 2), and $x+nm'\in W$ if and only if $x+nn'\in W$ by Condition 3).\\
Condition 2). We get $x+mm'y\in W$ if and only if $x+nm'y\in W$ by Condition 2), and we get $x+nm'y\in W$ if and only if $x+nn'y\in W$ by Condition 4).
Condition 3) is obtained in the same way from Conditions 3) and 4).\\
Condition 4). We get $x+ymm'z\in W$ if and only if $x+ynm'z\in W$ by Condition 4); and  $x+ynm'z\in W$ if and only if $x+ynn'z\in W$ again by Condition 4).
\endproof
\begin{theo}
	The congruence $R_W$ on $A$ is such that:\\
	1) the inclusion $w: W\into A$ is saturated w.r. to $R_W$;\\
	2) it is the largest congruence $S$ in $\SRg$ satisfying the saturation property w.r. to $W$. 
\end{theo}
\proof
1) Suppose $mR_Wn$. Then setting $x=0$ in Condition 1), we get $m\in W\iff n\in W$.\\
2) Suppose $S$ is a congruence w.r. to which $W$ is saturated, namely such that, when $mSn$ we have $m\in W \iff n\in W$. Now, from $mSn$, $\forall (x,y)\in X\times X$, we get $(x+m)S(x+n)$, $xmSxn$, $mySny$ and $xmySxny$. From that, $mR_Wn$ follows immediately. 
\endproof
We can now proceed exactly in the same way as in $\Mon$:
\begin{lemma}
	When we have $0\in W$, we then get: $\overline W\subset W$.
	So, when $W$ is a subsemiring, we get $\overline W\subset W$. 
\end{lemma}
Clearly, we get $W=\overline W$ under the following conditions: 
\begin{defi}
	Let $W$ be a subsemiring of $A$. Call it normal when for any $k\in W$, the following conditions  hold:\\
	1) $\forall x\in A,\; [x+k\in W \iff x\in W]; \;\;
	2) \; \forall x\in A, \; kx \in W \; {\rm and} \; xk\in W$.
\end{defi}

\begin{prop}
	A subsemiring $W$ of $A$ is normal if and only if it is a kernel inclusion in $\SRg$.
	
	Let $W$ be a normal subsemiring of $A$. Then $R_W$ is the supremum of the set of congruences $R$ on $A$ in $\SRg$ to which $W$ is normal.
\end{prop}

\subsubsection{Rings}

\begin{lemma}
	Let $R$ be a ring and $W\hookrightarrow R$ a subsemiring. Then $W$ is a normal subobject in $\SRg$ if and only if $W$ is an ideal of $R$.
\end{lemma}
\proof
Suppose $W$ is an ideal. Then Condition 2 of a normal subseming is satisfied by definition of an ideal, and Condition 1 holds since $(W,+)$ is a subgroup of $(A,+)$.
Conversely suppose $W$ a normal subsemiring and $k\in W$. Then $-k+k=0\in W \iff -k\in W$. And $(W,+)$ is a subgroup of $(R,+)$. Now, by the condition 2) of a normal subobject, $W$ is an ideal of $R$.
\endproof
Let $\Rng$ denote the category of rings, and $W$ be an ideal of a ring $R$. We know that the set of internal equivalence relations in $\Rng$ to which the ideal  $W$ is normal is reduced to one element; accordingly the category $\Rng$ is another extremal example of a category in which the set of internal equivalence relations to which a subobject is normal, when it is non-empty, has a supremum. Let us check that, in $\Rng$, this equivalence relation coincides with the syntactic equivalence relation $R_W$. Again, this is a corollary to the following
\begin{lemma}
	Let $R$ be a ring. Any equivalence relation $S$ on $R$ internal in $SRg$ is internal in $\Rng$.
\end{lemma}

\subsubsection{Left skew braces}

Braces are algebraic structures which were introduced by Rump \cite{Rump} as producing set-theoretical solutions of the Yang-Baxter equation in response to a general incitement of Drinfeld \cite{Drinfeld} to investigate this equation from a set-theoretical perspective. Later on, Guarnieri and Vendramin generalized this notion with the structure of left skew brace \cite{GV} which again generates solutions of the Yang-Baxter equations:\\
\noindent\textbf{Definition} 
\emph{A left skew brace is a set $X$ endowed with two group structures $(X,*,\circ)$ subject to a unique axiom: $a\circ (b*c)=(a\circ b)*a^{-*}*(a\circ c)\; (1)$, where $a^{-*}$ denotes the inverse for the law $*$.}

The simplest examples are the following ones: starting with any group $(G,*)$, then $(G,*,*)$ and $(G,*,*^{op})$ are left skew braces. We denote by $\SkB$ the category of left skew braces which is obviously a variety in the sense of Universal Algebra. Accordingly, it allows the definition of the syntactic equivalence relation, even though, in this case, the list of its axioms does not seem to be reductible to a finite one.

\section{Generalized syntactic equivalence relations}

Let us introduce the following definition whose some examples are briefly investigated at the end of this section. Section 5 will be devoted to produce some other ones.
\begin{defi}\label{supnorm}
	A category $\EE$ is said to be $\Eq$-saturating when, given any pair $(u,S)$ of a monomorphism $u: U\into X$ and an internal equivalence relation $S$ on $U$, the set of  saturated monomorphisms above $u$ with a domain smaller than $S$ has a supremum. We shall denote the codomain of this supremum by $\forall_uS$.
\end{defi}
From the universal property of $\forall_u$, we get immediately:
\begin{lemma}\label{fallsub}
	In any $\Eq$-saturating category $\EE$, we get:\\
1) for any monomorphism $u: U\into X$, then $T\subset S$ implies $\forall_uT\subset \forall_uS$;\\
2) when $v: X\into X'$ is another monomorphism, we get $\forall_v(\forall_uS)\subset \forall_{v.u}S$;\\
3) given any pullback:
\[\xymatrix@C=2pc@R=1pc{ {U\;} \ar[d]_{\bar f}  \ar@{>->}[r]^{f^{-1}(v)} & X \ar[d]^{f}\\
	{V\;} \ar@{>->}[r]_{v} & Y}\]
and any equivalence relation $S$ on $V$, we get: $f^{-1}(\forall_vS)\subset \forall_{f^{-1}(v)}(\bar f^{-1}S)$.
\end{lemma}
\begin{prop}
	Let $\EE$ be a $\Eq$-saturating category and $u: U\into X$ any monomorphism. Then an equivalence $R$ on $U$ is the domain of a saturated subobject above $u$ if and only if $R=u^{-1}(\forall_uR)$. Then $\forall_uR$ is the largest equivalence relation on $X$ saturated w.r. $R$ above $u$.
\end{prop}
\proof
Clearly when $R=u^{-1}(\forall_uR)$, $R$ is the domain of a saturated subobject. Conversely, suppose $u: R\into W$ is a saturated subobject. The universal property of $\forall_uR$ produces an inclusion $W\subset\forall_uR$, which makes $R=u^{-1}(\forall_uR)$. From that, the last assertion is straightforward.
\endproof
\begin{coro}\label{supnorm}
	Let $\EE$ be a $\Eq$-saturating category and $u: U\into X$ a  monomorphism. If $u$ is normal to some equivalence relation $T$ on $X$, the set of  the equivalence relations $S$ w.r. to which $u$ is normal has a supremum
\end{coro}
\proof 
Take $R=\nabla_U$ in the previous proof.
\endproof
\begin{theo}\label{var}
	Let $\EE$ be a $\Eq$-saturating category and $u: R\into S$ any monomorphism in $\Eq\EE$ above $u: U \into X$. Then $\overline S=S\cap (\forall_uR)$ is the largest of the equivalence relations $T$ on $X$ such that $u$ is saturated w.r. to $T$, $u^{-1}(T)\subset R$ and $T\subset S$.
\end{theo}
\proof
Since $\overline S\subset (\forall_uR)$, the monomorphism $u$ is saturated w.r. to $\overline S$ by Lemma \ref{lemma0}. And we have: $u^{-1}(\overline S)\subset u^{-1}(\forall_uR)\subset R$. Suppose we have an equivalence $T$ on $X$ satisfying the asserted properties. The two first ones implies $T\subset \forall_uR$, whence $T\subset \overline S$ with the third one. 
\endproof
The $\Eq$-saturating axiom has very good stability properties:
\begin{prop}\label{fibers}
Let $\EE$ be a $\Eq$-saturating category. Then any slice category $\EE/Y$, any coslice category $Y/\EE$ is $\Eq$-saturating as well. Accordingly any fiber $\Pt_Y\EE$ of split epimorphisms above $Y$ and any fiber $RGh_Y\EE$ of reflexive graphs on the object $Y$ is $\Eq$-saturating.
\end{prop}
\proof
For the coslice categories, it is straightforward since the left exact forgetful functor $cod: Y/\EE\to \EE$ determines a bijection between the equivalence relations on an object $h\in Y/\EE$ and the equivalence relations on its image by $cod$.

Let $f: X\to Y$ be an object of $\EE/Y$. An equivalence relation on this object in $\EE/Y$ is nothing but an equivalence relation $R$ on $X$ such that $R\subset R[f]$. Now consider any monomorphism $u$ in $\EE/Y$:
$$
\xymatrix@=15pt
{ {U\;} \ar[dr]_{f_u}\ar@{ >->}[rr]^{u}
	&& X\ar[dl]^{f} \\
	&  Y
}
$$
and any equivalence relation $S$ on the object $f_u$. Let us check that the equivalence relation $R[f]\cap (\forall_US)$ satisfies the desired universal property. Since $R[f]\cap (\forall_US)\subset R[f]$, this equivalence relation is on the object $f\in\EE/Y$. Now, the universal property follows from Theorem \ref{var}.

The last assertion comes from: $\Pt_Y\EE=1_Y/(\EE/Y)$  and from \\ $RGh_Y\EE=s_0^Y/(\EE/(Y\times Y))$, where $s_0^Y:Y\into Y\times Y$ is the diagonal.
\endproof                                                                 

\subsection{Monoids}

\begin{defi}
	Call {\em generalized syntactic relation} associated with the pair $(U,S)$ of a submonoid $U\hookrightarrow M$ and an internal equivalence relation $S$ on $U$, the relation $\forall_US$ on $M$ defined by: $$m(\forall_US)n \iff  \forall (x,y)\in M\times M,\;  [xmy\in U \iff xny\in U] \; {\rm and} \; xmySxny $$
\end{defi}
\begin{prop}\label{prop2}
	Given  any pair $(U,S)$, the relation $\forall_US$ is an internal equivalence relation in $\Mon$ whose class of the unit $1$ is the monoid $\overline U$ defined by:
	$$\overline U=\{ m\in M/ \forall (x,y)\in M\times M,\; [xmy\in U \iff xy\in U] \; {\rm and} \; xmySxy \}$$
\end{prop}
\proof
Same proof, modulo S, as in Proposition \ref{prop1}.
\endproof
\begin{theo}
	The category $\Mon$ is $\Eq$-saturating. 
\end{theo}
\proof
Given any pair $(U,S)$, we are going to show that the internal equivalence relation $\forall_US$ on $M$ satisfies the desired universal property:\\
1) the inclusion $u: U\into M$ is saturated w.r. to it and $u^{-1}(\forall_US)\subset S$:\\
When $(m,n)$ is in $U\times U$ and $m(\forall_US)n$, with $(x,y)=(1,1)$, we get $mSn$. Moreover when $m(\forall_US)n$ and $m\in U$, with the same $(1,1)$, we get $n\in U$.\\
2) $\forall_US$ is the largest internal equivalence relation on $M$ in $\Mon$ satisfying this saturation property:
suppose $u$ saturated w.r. to the equivalence relation $\Sigma$ on $M$ and $u^{-1}(\Sigma)\subset S$. Start with $m\Sigma n$. For all $(x,y)\in M\times M$, we get $xmy\Sigma xny$, and since $u$ is saturated w.r. to  $\Sigma$, we get $xmy\in U\iff xny\in U$. Now we get $xmySxny$ since $u^{-1}(\Sigma)\subset S$.
\endproof

\subsection{Groups}

In $\Gp$ it is easy to check that, given any pair $(u,R)$ of a monomorphism $u: U\into G$ and an equivalence relation $R$ on $U$, there is atmost one saturated subobject above $u$ with codomain $R$. We easily get a characterization:

\begin{prop}
	A monomorphism $u: U\into G$ in $\Gp$ is saturated w.r. to some equivalence relation $S$ on $G$ if and only if there exists a normal monomorphism (or normal subgroup) $v: V\into U$ such that the monomorphism $u.v$ is normal.
\end{prop}
\proof
Suppose $u: R\into S$ saturated in $\Gp$. Denote by $\iota: \bar 1^{R}\into U$ the normal subgroup produced by the class of the unit modulo $R$. Then it is  the class of the unit modulo $S$ in $G$. 

Conversely, suppose $v$ and $u.v$ are normal subobjects. Denote $R$ and $S$ the unique internal equivalence relations on $U$ and $X$ in $\Gp$ to which $v$ and $u.v$ are normal. Since $u.v$ is normal to $S$, the monomorphism $v$ is normal to $u^{-1}(S)$. Accordingly $R=u^{-1}(S)$, and $u: R\into S$ is saturated.
\endproof

Given any pair $(u, R)$ of a monomorphism $u: U\into G$ in $Gp$ and an equivalence $R$ on $U$ is $\Mon$, the {\em generalized syntactic equivalence relation $\forall_uR$ on $G$ lies in $Gp$} according to Lemma \ref{lemma2}.

\subsection{Semirings, left skew braces}

It is clear that the definition of the generalized syntactic relation with its characteristic universal property can be extended to any variety $\VV$.
\begin{defi}
	Given any algebra $A$ in $\VV$, any subalgebra $L$ and any congruence $S$ on $L$, call {\em generalized syntactic relation} associated with the pair $(L,S)$ the relation $R_L$ on $A$ defined by $mR_Ln$ when:\\ for any $(a_1,\cdots, a_n)\in A^n$ and any term $\tau (x_0,x_1,\cdots, x_n)$ of $\VV$, we get:\\
	$[\tau (m,a_1,\cdots, a_n)\in L \iff \tau (n,a_1,\cdots, a_n)\in L]$\\
	and $\tau (m,a_1,\cdots, a_n)S\tau (n,a_1,\cdots, a_n)$.
\end{defi}
From that, we get:
\begin{theo}
	When $L$ is a subalgebra of $A$, the relation $R_L$ is a congruence on $A$ in $\VV$ such that:\\
	1) the inclusion $l: L\into A$ is saturated w.r. to $R_L$;\\
	2) it is the largest congruence $T$ on $A$ in $\VV$ w.r. to which $l$ is saturated and such that $l^{-1}T\subset S$.\\ 
Accordingly, any variety $\VV$ is $\Eq$-saturating.
\end{theo}
\proof
Same proof, modulo $S$, as in Theorem \ref{alme}. Namely:\\
1) The relation $R_L$ is an equivalence relation since so are the hyperrelation [$\iff$] and the relation $S$.\\
2) To show that it is a congruence, start with a term  $\sigma(y_0,y_1,\cdots,y_m)$ and $mR_Ln$. We have to show that for any term $\tau (x_0,x_1,\cdots, x_n)$ we have: (i)\\
$\sigma(\tau (m,x_1,\cdots, x_n),y_1,\cdots,y_m)\in L\iff\sigma(\tau (n,x_1,\cdots, x_n),y_1,\cdots,y_m)\in L$
and (ii) $\sigma(\tau (m,x_1,\cdots, x_n),y_1,\cdots,y_m)S\sigma(\tau (n,x_1,\cdots, x_n),y_1,\cdots,y_m)$. Which is true since $\sigma(\tau (x_0,x_1,\cdots, x_n),y_1,\cdots,y_m)$ is a term and $mR_Ln$.\\
3) The inclusion $l: L\into A$ is saturated w.r. to $R_L$: use the term $\tau(x)=x$.\\
4) Universal property. Let $T$ be any congruence on $A$ such that $l$ is saturated w.r. to it and $l^{-1}(T)\subset S$. Starting with $mTn$ and any term $\tau (x_0,x_1,\cdots, x_n)$: we get $\tau (m,x_1,\cdots, x_n)T\tau (n,x_1,\cdots, x_n)$ and since $l$ is saturated w.r. to $T$, we get: $[\tau (m,x_1,\cdots, x_n)\in L \iff \tau (n,x_1,\cdots, X_n)\in L]$, namely:\\ $[\tau (m,x_1,\cdots, x_n)l^{-1}(T)\tau (n,x_1,\cdots, x_n)]$. Then $mR_Ln$, thanks to $l^{-1}T\subset S$.
\endproof
So, the categories $\SRg$ of semi-rings, $\Rng$ of rings and $\SkB$ of left skew braces are examples of $\Eq$ saturating categories.

\section{The case of protomodular categories}

In this section we are going to investigate the properties of the $\Eq$-saturating protomodular categories as is the category $\Gp$. Recall that a (finitely complete category) $\EE$ is protomodular \cite{Bprot} when, given any pair of commutative squares of vertical split epimorphisms: 
$$
\xymatrix@=20pt{
	\bullet \ar[d] \ar[r] & \bullet \ar[d] \ar[r]  & \bullet \ar[d] \\
	\bullet \ar@<-1ex>@{ >->}[u]  \ar[r] & \bullet \ar[r] \ar@<-1ex>@{ >->}[u] & \bullet \ar@<-1ex>@{ >->}[u] 
}
$$
the right hand side square is a pullback as soon as so are the left hand side one and the whole rectangle. Denote by $\Pt\EE$ the category whose objects are the split epimorphisms $(f,s):X\splito Y$ and whose morphisms are the commutative squares between them. Let $\P_\EE: \Pt\EE \to \EE$ be the functor associating with any split epimorphism $(f,s)$ its codomain $Y$; it is a fibration (called {\em fibration of points}) whose cartesian maps are the pullbacks; we denote by $\Pt_Y\EE$ the fiber of split epimorphisms above $Y$. A category $\EE$ is protomodular if and only if any base-change functor of this fibration is conservative.

The major examples are the categories $\Gp$ of groups, $\Rng$ of rings, and any category of $R$-algebras when $R$ is a ring. The category $\SkB$ is protomodular as well, see \cite{BFP}. The protomodular varieties are characterized in \cite{BJ}. In a pointed protomodular category we get any of the classical homological lemmas, see \cite{BB}.
 
Let $\EE$ be a protomodular category. Then we get:\\
1) the class of fibrant morphisms in $\Eq\EE$  is such that, when $g.f$ and $f$ are fibrant, so is $g$ (it is an immediate consequence of the definition);\\
2) any fibrant morphism $f: R\to S$ in $\Eq\EE$ is cocartesian w.r. to the fibration $(\;)_0:\Eq\EE \to \EE$ \cite{BM};\\
3) a monomorphism $u: U\into X$ is normal to at most one equivalence relation $R$ on $X$ \cite{Bnorm};\\
4) an object $X$ is called {\em affine} when the diagonal $s_0^X: X\into X\times X$ is normal; this determines the protomodular subcategory $\Aff\EE\into\EE$ of affine objects.
\begin{lemma}\label{prot}
	Let $\EE$ be any protomodular category.\\
	1) Let $u: W\into R$ and  $u: W'\into R'$ be two saturated monomorphisms. Then we get: $[R'\subset R\iff W'\subset W]$.\\
	2) Let $u: W\into R$ be a saturated monomorphism. Then $R$ has the universal property of $\forall_UW$.\\
	3) There is (up to isomorphism) at most one saturated subobject $u: W \into R$ with domain $W$.
\end{lemma}
\proof
1) The direct implication is straightforward. Suppose $W'\subset W$. Then $u: W'\cap W=W' \into R'\cap R$ is saturated. The category $\EE$ being protomodular, we get $R'\simeq R'\cap R$ and $R'\subset R$. Then 2) follows immediately, and the associated universal property implies 3).
\endproof
\noindent The point 3) of this lemma obviously generalizes the observation 3) above it.

\subsection{$\Eq$-saturating protomodular categories}

The categories $\Gp$, $\Rng$, $\SkB$ produce examples of $\Eq$-saturating protomodular categories. This is more generally the case of any protomodular variety. We shall produce some other examples in the next sections.

\begin{lemma}
	Let $\EE$ be an  $\Eq$-saturating protomodular category and $(u,v)$ any pair of composable monomorphisms: $U\stackrel{u}{\into} X \stackrel{v}{\into} X'$. We then get $\forall_v(\forall_uS)=\forall_{v.u}S$.
\end{lemma}
\proof
Let $v.u: W\into T$ be a saturated subobject with $W\subset S$. Then $u: W\into v^{-1}T$ and $v: v^{-1}T\into T$ are saturated as well; $u: W\into v^{-1}T$ since $v.u: W\into T$ is saturated and $v: v^{-1}T\into T$ is a monomorphism, and then $v: v^{-1}T\into T$ by the observation 1) above the previous Lemma.

The monomorphism $u: W\into v^{-1}T$ being saturated such that $W\subset S$, we get an inclusison $v^{-1}T \subset \forall_uS$. From this inclusion and the saturated monomorphism $v: v^{-1}T\into T$, we get $T\subset \forall_v(\forall_uS)$. Whence $\forall_v(\forall_uS)=\forall_{v.u}S$.
\endproof
\begin{prop}
	Let $\EE$ be an  $\Eq$-saturating protomodular category.\\ 
	1) Let $u: U\into X$ be a monomorphism and $R$ any equivalence relation on $X$. We then get  $\forall_u(u^{-1}R)\subset R$ in such a way that the following square is a pullback in $\Eq\EE$:
	$$
	\xymatrix@=15pt{
		u^{-1}(\forall_u(u^{-1}R)) \ar@{ >->}[r] \ar@{ >->}[d]   & \forall_u(u^{-1}R) \ar@{ >->}[d]  \\
		u^{-1}R  \ar@{ >->}[r]  & R  
	}
	$$
	2) Let $u: S\into R$ be any monomorphism of equivalence relations. We then get: $(\forall_{u}S)\subset R$ in such a way that the following square is commutative in $\Eq\EE$:
		$$
	\xymatrix@=15pt{
		u^{-1}(\forall_uS) \ar@{ >->}[r] \ar@{ >->}[d]   & \forall_uS \ar@{ >->}[d]  \\
		S  \ar@{ >->}[r]  & R  
	}
	$$
\end{prop}
\proof
1) We get $u^{-1}(R\cap (\forall_u(u^{-1}R))=u^{-1}R\cap u^{-1}(\forall_u(u^{-1}R))=u^{-1}(\forall_u(u^{-1}R))$ since $u^{-1}(\forall_u(u^{-1}R)\subset u^{-1}R$. According to Lemma \ref{prot}, we get $R\cap (\forall_u(u^{-1}R))=\forall_u(u^{-1}R)$ and $\forall_u(u^{-1}R)\subset R$.\\
2) From $S\subset u^{-1}R$, we get $\forall_uS\subset \forall_u(u^{-1}R)$ by Lemma \ref{fallsub}; whence, by 1), $\forall_{u}S\subset R$. 
\endproof

\subsection{$\Eq$-saturating exact protomodular categories}\label{normali}

When moreover $\EE$ is an exact category \cite{Barr}, as is any protomodular variety, this result can be translated in terms of  a property of the fibration of points $\P_\EE$: 
\begin{prop}
	Let $\EE$ be an exact $\Eq$-saturating protomodular category. Any commutative square of vertical split epimorphisms with horizontal regular epimorphisms $(h,k)$ has a universal dotted decomposition:
	$$
	\xymatrix@=20pt{
		{X\;} \ar[d]_{f} \ar@<2ex>@{->>}[rr]^{k} \ar@{.>>}[r]_{\bar k} & {\bar X\;} \ar[d]_{\bar f} \ar@{.>>}[r]_{\check k}  & V \ar[d]_{g} \\
		{Y\;} \ar@<-1ex>[u]_{s} \ar@<-2ex>@{->>}[rr]_{h} \ar@{.>>}[r]^{\bar h} & {\bar Y\;} \ar@{.>>}[r]^{\check h} \ar@<-1ex>[u]_{\bar s} & W \ar@<-1ex>[u]_{t} 
	}
	$$
	where the left hand side square is a pullback with a regular epimorphism $\bar h$.
\end{prop} 
\proof
Consider the inclusion $s: R[h] \into R[k]$ and then the saturated inclusion $s: s^{-1}(\forall_sR[h]) \into (\forall_sR[h])$. Thank to the previous proposition we get $(\forall_sR[h])\subset R[k]$. According to the point 2) of the above presentation of protomodular categories, the fibrant map $s: s^{-1}(\forall_sR[h]) \into (\forall_sR[h])$ being cocartesian, the retraction $f$ of $s$ produces a map $f: (\forall_sR[h]) \to s^{-1}(\forall_sR[h])$ which is fibrant as well by the point 1).

As soon as $\EE$ is exact, take $\bar k$ and $\bar h$ the quotient maps of $s^{-1}(\forall_sR[h])$ and $(\forall_sR[h])$. Let $\bar s$ and $\bar f$ the induced factorizations through the quotients. By the Barr-Kock Theorem \cite{BK} valid in any exact category, the fact that the map $f: (\forall_sR[h]) \to s^{-1}(\forall_sR[h])$ is fibrant makes the left hand side downward square a pullback, and the upward as well. The inclusions $(\forall_sR[h])\subset R[k]$ and $s^{-1}(\forall_sR[h])\subset R[h]$ produce the desired factorizations $\check k$ and $\check h$.

Suppose now there is another decomposition:
	$$
\xymatrix@=20pt{
	{X\;} \ar[d]_{f} \ar@<2ex>@{->>}[rr]^{k} \ar@{.>>}[r]_{\bar\kappa} & {X'\;} \ar[d]_{f'} \ar@{.>>}[r]_{\check\kappa}  & V \ar[d]_{g} \\
	{Y\;} \ar@<-1ex>[u]_{s} \ar@<-2ex>@{->>}[rr]_{h} \ar@{.>>}[r]^{\bar\chi} & {Y'\;} \ar@{.>>}[r]^{\check\chi} \ar@<-1ex>[u]_{s'} & W \ar@<-1ex>[u]_{t} 
}
$$
where the left hand side square is a pullback. Then the inclusion $s: R[\bar\chi]\into R[\bar\kappa]$ is saturated. The factorization $\check\chi$ implies $R[\bar\chi]\subset R[h]$ which determines an inclusion $R[\bar\kappa] \subset (\forall_sR[h])$ and the desired factorization $\tilde\kappa: X'\to \bar X$, while $R[\bar\chi]=s^{-1}(R[\bar\kappa])\subset s^{-1}(\forall_sR[h])=R[\bar h]$ produces the desired factorization $\tilde\chi: Y' \to \bar Y$.
\endproof

This decomposition appears to be complementary to the one which characterises {\em the existence of normalizers} in quasi-pointed protomodular categories \cite{Bnormlizer}, namely the existence of a universal decomposition for the monomorphisms between split epimorphisms in $\mathbb E$: any monomorphism $(y,x):(f's')\rightarrowtail (f,s)$ of split epimorphisms produces a universal dotted decomposition where the left hand side part is a pullback:
$$
\xymatrix@=20pt{
	{X'\;} \ar[d]_{f'} \ar@<2ex>@{>->}[rr]^{x} \ar@{>.>}[r]_{\bar x} & {\bar X\;} \ar[d]_{\bar f} \ar@{>.>}[r]_{\check x}  & X \ar[d]_{f} \\
	{Y'\;} \ar@<-1ex>[u]_{s'} \ar@<-2ex>@{>->}[rr]_{y} \ar@{>.>}[r]^{\bar y} & {\bar Y\;} \ar@{>.>}[r]^{\check y} \ar@<-1ex>[u]_{\bar s} & Y \ar@<-1ex>[u]_{s} 
}
$$

\subsection{Connected pairs of equivalence relations}

Recall from \cite{BG} that, in a Mal'tsev category (in which any reflexive relation is an equivalence relation \cite{CLP}, \cite{CPP}), two equivalence relations  $(R,S)$ on an object $X$  {\em are connected} or {\em centralize each other} when there is a (necessarily unique) double equivalence relation on then:
$$ \xymatrix@=30pt{
	W \ar@<-1,ex>[d]_{p_0^R}\ar@<+1,ex>[d]^{p_1^R} \ar@<-1,ex>[r]_{p_0^S}\ar@<+1,ex>[r]^{p_1^S}
	& S \ar@<-1,ex>[d]_{d_0^S}\ar@<+1,ex>[d]^{d_1^S} \ar[l]\\
	R \ar@<-1,ex>[r]_{d_0^R} \ar@<+1,ex>[r]^{d_1^R} \ar[u]_{} & X
	\ar[l]|{s_0^R} \ar[u]
}
$$
such that the morphism $d_0^R: W\to S$ is fibrant. Then any commutative square in this diagram is a pullback. We classically denote this property by $[R,S]=0$. Recall that any protomodular category is a Mal'tsev one \cite{BB}.

\begin{prop}
	Let $\EE$ be a protomodular category, and $(R,S)$ a pair of equivalence relations on $X$. Then $[R,S]=0$ if and only if $S$ is the domain of a saturated monomorphism above $s_0^R$.
\end{prop}
\proof
Suppose $[R,S]=0$. Then, any square in the here above diagram being a pullback, the morphism $s_0^R: S\into W$ is fibrant. Conversely suppose $S$ is the domain of a saturated monomorphism $S\into W$ above $s_0^R$:
$$ \xymatrix@=30pt{
	W \ar@<-1,ex>[d]_{p_0^R}\ar@<+1,ex>[d]^{p_1^R} \ar@<-1,ex>@{.>}[r]_{p_0^S}\ar@<+1,ex>@{.>}[r]^{p_1^S}
	& S \ar@<-1,ex>[d]_{d_0^S}\ar@<+1,ex>[d]^{d_1^S} \ar[l]\\
	R \ar@<-1,ex>@{.>}[r]_{d_0^R} \ar@<+1,ex>@{.>}[r]^{d_1^R} \ar[u]_{} & X
	\ar[u]_{} \ar[l]|{s_0^R}
}
$$
Since the morphism $s_0^R$ is cocartesian, it induces the morphism $d_i^R: W\to S$ which are fibrant, since the fibrant morphisms are right cancellable.
\endproof

\begin{theo}
	In any $\Eq$-saturating protomodular category $\EE$, any equivalence relation $R$ has a centralizer.
\end{theo} 
\proof
Let $R$ be an equivalence relation on the object $X$, then consider the equivalence relation $\forall_{s_0^R}\nabla_X$ on $R$ and the induced saturated monomorphism above $s_0^R$:
$$ \xymatrix@=20pt{
    \forall_{s_0^R}\nabla_X \ar@<-1,ex>[d]_{p_0^R}\ar@<+1,ex>[d]^{p_1^R} 
	& \Sigma \ar@<-1,ex>[d]_{d_0^\Sigma}\ar@<+1,ex>[d]^{d_1^\Sigma} \ar[l]\\
	R  \ar[u]_{} & X
	\ar[u]_{} \ar[l]|{s_0^R}
}
$$
According to the previous proposition, we get $[R,\Sigma]=0$. Now, if we have $[R,S]=0$ for some $S$ with some $W$ as in the previous proposition, the universal property of $\forall_{s_0^R}\nabla_X$ implies the inclusion $W\subset \forall_{s_0^R}\nabla_X$ and thus $S\subset \Sigma$.
\endproof
\noindent The existence of centralizers in the category $\SkB$ was shown in \cite{BGG}, but in a more classical way. Exactly in the same way, we get:
\begin{prop}
	Let $\EE$ be an  $\Eq$-saturating protomodular category and $R$ any equivalence relation on an object $X$. Given any other equivalence relation $S$ on $X$, then $(s_0^R)^{-1}(\forall_{s_0^R}S)$ is the largest of the equivalence relations $T$ on $X$, such that $[R,T]=0$ and $T\subset S$.
\end{prop}

\subsection{The pointed protomodular and additive cases}

In the pointed context, we get a preorder homomorphism $\iota_X:\Eq_X\EE \to \Sub X$, where $\Sub X$ is the set of monomorphisms with codomain $X$. Starting from an equivalence $R$ on $X$ we associate with it the following the following upper horizontal monomorphism, where the left hand side square is a pullback:
$$ \xymatrix@=20pt{
 I_R \ar[d]_{} \ar@{ >->}[r]_{\alpha_R} \ar@<2ex>@{ >->}[rr]^{\iota_X^R} &	R \ar[d]^{d_0^R}\ar[r]_{d_1^R} 	& X \\
 1 \ar@{ >->}[r]_{\alpha_X} &	X  
}
$$
\begin{lemma}
	Given any pointed category $\EE$ and any saturated monomorphism $u: R\into S$, we get $I_U^R\simeq_X^S$.
\end{lemma}
\proof
Straightforward from the stability of pullbacks by composition.
\endproof
As soon as the pointed category $\EE$ is protomodular, the   preorder homomorphism $\iota_X$ is full (i.e. reflects inclusions) by Lemma \ref{prot}.1, and its image determines the preorder $\Norm\EE$ of normal monomorphisms with codomain $X$.
We get moreover:
\begin{lemma}
	Let $\EE$ be a pointed protomodular category and $u: R\into S$ a monomorphism in $\Eq\EE$. It is saturated if and only if $I_U^R\simeq I_X^S$. In other words, a monomorphism $u: U\into X$ is saturated w.r. to $S$ if and only if $I_X^S\subset U$.
\end{lemma}
\proof
The condition is necessary by the previous lemma. 
Conversely, the monomorphism $u: R\into S$ induces an inclusion $j: I_U^R\into I_X^S$. More precisely, it determines the following quadrangled  pullback of split epimorphisms since $\tilde u$ is a monomorphism:
$$\xymatrix@=7pt{
	{I_U^R\;} \ar@{>->}[rd]^{j} 
		  \ar@<-1ex>[ddd]_<<<<{}  \ar@{ >->}[rrrr]^{\alpha_R}  &&&&  {R\;} \ar@{ >->}[rrrrrd]^{\tilde u}  \ar@<-1ex>[ddd]_>>>>>>>>{d_0^R} \ar@{ >.>}[rd] &&&&   \\
	&   I_X^S \ar[ddl]_{} \ar@{ >->}[rrrr]^<<<<{}  && && \bullet \ar@{ >->}[rrrr]^<<<<{}\ar@{.>}[ddl]_{}  && && S \ar[ddl]_{d_0^S} \\
	&&& &&&&&&\\
	1  \ar@{ >->}[rrrr]_{\alpha_U}  \ar[uuu]_>>>>{} \ar@<-1ex>[uur]_{} && && U \ar@{ >->}[rrrr]_{u} \ar@<-1ex>@{.>}[uur]_{} \ar[uuu]_<<<<<<<<{s_0^R} && && X \ar@<-1ex>[uur]_{s_0^S} 
}
$$
If $\bullet$ denotes the pullback of $d_0^S$ along $u$, then the induced dotted factorization $\bar j:R\into \bullet$ produces a parallelistic pullback. Since $\EE$ is protomodular and $j$ an isomorphism, so is $\bar j$, and $u: R\into S$ is saturated.
\endproof
So, in the pointed protomodular context, this translation allows a direct access to the $\Eq$-saturating property:
\begin{prop}
Let $\EE$ be a pointed protomodular category. It is $\Eq$-saturating if and only if, given any monomorphism $u: U\into X$, any normal subobject $v: V\into U$ determines a largest normal subobject $w: W\into X$  such that $W\subset V$.	
\end{prop}

\medskip

The quickest way to define an {\em additive category} is to say that any object $X$ is endowed with a natural internal group structure we shall denote $(X,+)$ since, being natural, this group structure is then necessarily unique and abelian. So, any additive category is pointed. For any pair $(X,Y)$ of objects, we shall denote by $0: X\to Y$, the map $X\stackrel{\tau_X}{\to}1 \stackrel{\alpha_Y}{\to} Y$. Recall from \cite{Bnorm}:
\begin{prop}
	Let $\EE$ be a pointed protomodular category.\\
	1) an object $X$ is abelian if and only if the diagonal $s_0^X: X \into X\times X$ is normal;\\ 
	2) $\EE$ is additive if and only if any diagonal is normal; in this case, any monomorphism is normal.
\end{prop}
\begin{prop}
	Given any additive category $\EE$, the full preorder homomorphism  $\iota_X:\Eq_X\EE \into \Sub X$ is an equivalence.
\end{prop}
\proof
With any suboject $v: V\into X$, you can associate the following reflexive relation $T_v$:
$$
\xymatrix@=20pt
{X\times V \ar@<+2ex>[rr]^-{+.(X\times v)}\ar@<-2ex>[rr]_-{p_X}
	&& X\ar[ll]|<<<<<<<<{(1_X,0)}
}$$
which determines an equivalence relation on $X$ since any protomodular category is a Mal'tsev one.
\endproof
Accordingly, in the additive context, the interest of the $\Eq$-saturating property vanishes since it becomes trivial.
\begin{prop}\label{addi}
	Any additive category $\EE$is (trivially) $\Eq$-saturating.
\end{prop}

Thanks to the previous equivalence, given any monomorphism $u: U\into X$, we have to find, for any monomorphism $v: V\into U$, a largest monomorphism $W: W\into X$ such that $W\subset V$. Clearly this largest monomorphism is $u.v: V\into X$.
\endproof

\section{The fibers of the fibration $(\;)_0: \Cat\EE \to \EE$}

The aim of this section is to produce new examples of $\Eq$-saturating categories.\\
Given any reflexive graph $X_{1}$:
$$
\xymatrix@=20pt
{X_1\ar@<+2ex>[r]^-{d_1}\ar@<-2ex>[r]_-{d_0}
	& X_0\ar[l]|{s_0}
}
$$
in a category $\EE$, it is underlying an internal category $X_{\bullet}$ in $\EE$ when there is moreover a composition map $d_1^{X_{\bullet}}: X_1\times_0X_1 \to X_1$ on the composable pairs in $X_{1}$ which is associative and makes neutral the composition with the identity maps, see for instance \cite{Bcat} for more details. Whence the category $\Cat\EE$ of internal categories and internal functors, and a left exact functor: $(\;)_0: \Cat\EE \to \EE$ associating with any internal  category $X_{\bullet}$ its "object of objects" $X_0$. It is a fibration whose cartesian maps are the fully faithful internal functors. The fiber above $X$ is denoted by $\Cat_X\EE$. An internal groupoid is an internal category $X_{\bullet}$ such that "any morphism is invertible", namely such that the following commutative square is a pullback in $\EE$:
$$
\xymatrix@=20pt
{X_1\times_0X_1\ar[r]^{d^{X_{\bullet}}_0} \ar[d]_{d^{X_{\bullet}}_1}
	& X_1 \ar[d]^{d^{X_{\bullet}}_0}\\
	X_1\ar[r]_-{d^{X_{\bullet}}_0}
	& X_0
	}
$$
The internal groupoids are  stable under cartesian maps and determine a subfibration $(\;)_0:\Grd\EE \to \EE$ of the previous one. Any fiber $\Grd_X\EE$ is a protomodular category, see \cite{Bprot}, as is $\Grd_1\EE=\Gp\EE$ the category of internal groups in $\EE$.

We are going now to produce examples of category $\EE$ with $\Eq$-saturating fibers $\Cat_X\EE$ and $\Grd_X\EE$.

\subsection{The fibration $(\;)_0: \Cat \to \Set$}

Let $\Cat$ be the category of categories, $\Grd$ the subcategory of groupoids and  $(\;)_0:\Cat \to \Set$ the associated fibration whose cartesian maps are the fully faithful functors. We denote by $\Cat_X$ and $\Grd_X$ the respective fibers. The respective fibers above the singleton $1$ are $\Mon$ and $\Gp$. Similarly to $\Gp$, any fiber $\Grd_X$ is protomodular. 

Let $\gamma: \CC \into \DD$ be a bijective on objects inclusion ($C_0=X=D_0$), and $S$ an internal equivalence relation on $\CC$  in the fiber $\Cat_X$, namely an equivalence relation on parallel pairs of morphims.
Define the following {\em generalized syntactic} relation $\forall_\CC S$ on the parallel pairs $(f,f'): a\rightrightarrows b$ of morphisms in $\DD$: $f(\forall_\CC S)f'$ when
$\forall (g,k)$ pair of maps in $\DD$ with $dom(g)=b$ and $cod(k)=a$, we get: $$[gfk\in \CC \iff gf'k\in \CC] \; {\rm and} \; gfkSgf'k $$

\begin{prop}
	Given any set $X$, the fiber $\Cat_X$ is an $\Eq$-saturating categories.
\end{prop}
\proof
Let us check that the equivalence relation $\forall_\CC S$ on $\DD$ is internal in $\Cat_X$ and that it has the desired universal property. It is clearly an equivalence relation on the pair of parallel maps in $\DD$. Let us check it is internal in the fiber $\Cat_X$. Suppose $f(\forall_\CC S)f'$ and $h(\forall_\CC S)h'$ with $dom(h)(=dom(h'))=b$ and $cod(h)(=cod(h'))=c$. We have to check $h.f(\forall_\CC S)h'.f'$. As in Proposition \ref{prop1}, $\forall (g,k)$ with $cod(k)=a$ and $dom(g)=c$, we get $[ghfk\in \CC \iff ghf'k \in \CC]$ with $ghfkSghf'k$, since $f(\forall_\CC S)f'$. And 
$[ghf'k\in \CC \iff gh'f'k \in \CC]$ with $ghf'kSgh'f'k$, since $h(\forall_\CC S)h'$. Whence $[ghfk\in \CC \iff gh'f'k \in \CC]$ and $ghfkSgh'f'k$.
It remains to check that $\forall_\CC S$ has the desired universal property.

First, let us check that the inclusion $\CC \into \DD$ is saturated w.r. to $\forall_\CC S$. Suppose $\phi(\forall_\CC S)f$ with $\phi\in \CC$. Taking the pair $(1_b,1_a)$ we get: $\phi\in \CC$ if and only if $f\in \CC$; and $\phi Sf$, which means $\gamma^{-1}(\forall_\CC S)\subset S$.

Suppose now $T$ is an internal equivalence relation in $\Cat_X$ on $\DD$ w.r. to which the inclusion $\gamma$ is saturated and is such that $\gamma^{-1}(T)\subset S$. Suppose $fTf'$ in $\DD$. We get $gfkTgf'k$ for all pair $(g,k)$ of composable maps in $\DD$, and, since the inclusion is saturated w.r. to $T$, we get $[gfk\in \CC \iff gf'k\in \CC]$, and $gfkSgf'k$ since $\gamma^{-1}(T)\subset S$. 
\endproof
\begin{coro}
	Given any set $X$, the fiber $\Grd_X$ is a $\Eq$-saturating protomodular category.
\end{coro}
\noindent It is a consequence of the following
\begin{lemma}\label{lemma3}
	Let $\DD$ be a groupoid and $S$ an equivalence relation on $\DD$ which is internal in the fiber $\Cat_X$. Then $S$ is internal in the fiber $\Grd_X$.
\end{lemma}
\proof
It remains to check that when we have $fSf'$ in $\DD$, we get $f^{-1}Sf'^{-1}$. From $fSf'$, we get $1_b=f.f^{-1}Sf'.f^{-1}$; then $f'^{-1}Sf^{-1}$.
\endproof

\subsection{The case of Mal'tsev and Gumm categories}

\subsubsection{Mal'tsev categories}\label{mal}

A Mal'tsev category is, as we recalled it, a category in which any reflexive  relation is an equivalence relation, see \cite{CLP} and \cite{CPP}. Equivalently it is a category in which any subgraph of an internal groupoid is necessarily a groupoid \cite{B1}. On the other hand, in a Mal'tsev category, any internal category is a groupoid \cite{CPP}, so that any fiber $\Cat_X\EE=\Grd_X\EE$ is protomodular.

Furthermore, any fiber $\Grd_X\EE$ is a naturally Mal'tsev category \cite{abgrd}, where a naturally Mal'tsev category \cite{Johnstone} is a category $\EE$ such that any object $X$ is endowed with a natural internal Mal'tsev operation, namely a ternary operation $\pi_X: X\times X\times X \to X$ such that $\pi_x(x,y,y)=x=\pi_X(z,z,x)$. 
In this way, a naturally Mal'tsev category is a first step inside a non-pointed additive context. Any naturally Mal'tsev category is a Mal'tsev one, and from \cite{BB}:
\begin{prop}
	Given any Mal'tsev category $\EE$, the following conditions are equivalent:\\
	1) it is a naturally Mal'tsev one;\\
	2) any equivalence relation $R$ on $X$ is central, namely $[R,\nabla_X]=0$;\\
	3) any object $X$ is affine, namely $[\nabla_X,\nabla_X]=0$.
\end{prop}

In \cite{abgrd} was introduced the following classification table of the different degrees of non-pointed additiveness by decreasing order of generality and with illustrating examples:
$$\begin{tabular}{|c|c|c|}
\hline
category $\EE$& fibration $\P_\EE :\Pt\EE \to \EE$ & examples \\
\hline
\hline
&  & $\Aff\EE$  \\
nat. Mal'tsev & additive fibers & when $\EE$ Mal'tsev\\
\hline
protomodular & conservative base-change & $\Aff\EE$ \\
+ nat. Mal'tsev &functors  + additive fibers &  when $\EE$ proto \\
\hline
antepenessentially & fully faithful & $\Grd_X\EE$\\
affine & base-change functors & when $\EE$ Gumm \\
\hline
 penessentially &  fully faithful + saturated & $\Grd_X\EE$ \\
  affine & on subobj. base-change functors  & when $\EE$ Mal'tsev\\
  \hline
 essentially affine & base-change functors are & $\Grd_X\EE$\\
  & equivalences of categories & when $\EE$ additive\\
\hline
\end{tabular}$$ 
The three last cases are protomodular. Generalizing Proposition \ref{addi} to the non-pointed case, we get:
\begin{prop}
	Any penessentially affine category  $\EE$ is $\Eq$-saturating. In such a category, any monomorphism is normal.
\end{prop} 
\proof
Let $u: U\into X$ be any monomorphism and $S$ any equivalence relation on the object $U$. Consider the following diagram with the monomorphism $j=(d_0^S,u.d_1^S)$:
$$\xymatrix@=9pt{
	{S\;} \ar@{>->}[rd]^{j} 
	\ar@<-1ex>[ddd]_<<<<<<{d_0^S}  \ar@{ >.>}[rrrr]^{\tilde u}  &&&&  {R\;}  \ar@<-1ex>[ddd]_>>>>>>{d_0^R} \ar@{ >.>}[rd]^{\bar j} &&&&   \\
	&   U\times X \ar[ddl]_{p_U} \ar@{ >->}[rrrr]_<<<<<<{u\times X}  && && X\times X \ar[ddl]_{p_0^X}   \\
	&&& &&&&&&\\
	U  \ar@{ >->}[rrrr]_{u}  \ar[uuu]_>>>>>>{s_0^S} \ar@<-1ex>[uur]_{(1_U,u)} && && X  \ar@<-1ex>[uur]_{s_0^X} \ar[uuu]_>>>>>>{} 
}
$$
Since the base-change $u^*$ is saturated on subobjects its determines a subobject $\bar j: (d_0^R,s_0^R)\into (p_0^X,s_0^X)$ in the fiber $Pt_X\EE$ of split epimorphisms above $X$ such that $u^*((d_0^R,s_0^R))=(d_0^S,s_0^S)$ ($*$) and $u^*(\bar j)=j$. We get the following reflexive relation which is an equivalence relation since any penessentially affine category is protomodular:
$$\xymatrix@=20pt{
R \ar@<+2ex>[r]^{p_1^X.\bar j}\ar@<-2ex>[r]_{d_0^R}
	& X \ar[l]|{s_0^R}
}
$$
With ($*$) and $u^*(d_0^R,s_0^R)=(d_0^S,s_0^S)$, the monomorphism $u: S\into R$ is saturated, and according to Lemma \ref{prot}.2, we get $R=\forall_uS$. For the second assertion, start with $S=\nabla_U$. This last point was already shown in \cite{abgrd}
\endproof
\begin{coro}
	When $\EE$ is a Mal'tsev category, any fiber $\Cat_X\EE=Grd_X\EE$ is an $\Eq$-saturating protomodular category.
\end{coro} 

\subsubsection{Gumm categories}

A Gumm category \cite{BGG} is a category in which the {\em Shifting Lemma} holds: {\em given any triple $(T,S,R)$ of equivalence relations on an object $X$ such that} $R\cap S\subset T$, {\em the following left hand side situation implies the right hand side one}:
$$ \xymatrix@=10pt{
	x \ar[r]^{S} \ar@(l,l)[d]_T \ar[d]_R & y \ar[d]^R &&  y \ar@(r,r)[d]^T\\
	x' \ar[r]_{S} & y' &&   y'
}
$$ 
The name comes from this characterization by Gumm of the congruence modular varieties \cite{Gu}.
Let us recall important results about internal categories in Gumm categories $\EE$ \cite{BGgumm}: on a reflexive graph $X_{\bullet}$, there is at most one structure of internal category, and the induced inclusion $\Cat\EE\hookrightarrow \RGh\EE$ is a full inclusion. Moreover, in \cite{Bgumm}, we get the following:
\begin{theo}
	Let be given an internal category $X_{\bullet}$ in a Gumm category $\EE$\linebreak
	$$\xymatrix@=6pt{
		&&& R_1\ar@<-4pt>[dd]_{d_0} \ar@<4pt>[dd]^{d_1} \ar@<-4pt>[rrr]_{d_0} \ar@<4pt>[rrr]^{d_1} &&& R_0 \ar@{ >->}[lll]|{s_0} \ar@<-4pt>[dd]_{d_0} \ar@<4pt>[dd]^{d_1} &&&&& R_{\bullet}\ar@<-4pt>[dd]_{d_0} \ar@<4pt>[dd]^{d_1}\\
		&&&&&&&&&&&\\
		X_1\times_0 X_1 \ar@<-4pt>[rrr]_{d_0} \ar[rrr]|{d_1} \ar@<4pt>[rrr]^{d_2}  &&&	X_1 \ar@{ >->}[uu]|{} \ar@<-4pt>[rrr]_{d_0} \ar@<4pt>[rrr]^{d_1} 
		&&& X_0 \ar@{ >->}[lll]|{s_0}    \ar@{ >->}[uu] &&&&&  X_{\bullet}  \ar@{ >->}[uu]
	}
	$$
	together with a vertical internal equivalence relation $R_{\bullet}$ on the underlying reflexive graph of $X_{\bullet}$.
	Then the upper horizontal reflexive graph $R_{\bullet}$ is underlying an internal category. Namely, $R_{\bullet}$ produces an  internal equivalence  relation in $\Cat\EE$.
\end{theo}
 In other words, the previous full inclusion $\Cat\EE\hookrightarrow \RGh\EE$ is stable under equivalence relations. It is a fortiori the case for any full inclusion: $\Cat_X\EE\hookrightarrow \RGh_X\EE$.
\begin{prop}
  When $\EE$ is an $\Eq$-saturated Gumm category, any fiber $\Cat_X\EE$ is an $\Eq$-saturating category as well.
\end{prop}
\proof
According to Proposition \ref{fibers}, any fiber $RGh_X\EE$ is $\Eq$-saturated. Since $\EE$ is a Gumm category, we get a full inclusion $\Cat_X\EE\hookrightarrow \RGh_X\EE$ which preserves and reflects the saturated subobjects between equivalence relations and is moreover stable under equivalence relations. The $\Eq$-saturating property of the fiber $\Cat_X\EE$ follows immediately.
\endproof
\begin{coro}
	Given an $\Eq$-saturating Gumm category $\EE$, any fiber $\Grd_X\EE$ is a $\Eq$-saturating protomodular category.
\end{coro}
\noindent Accordingly, these two results hold in any congruence modular variety $\VV$.

\medskip

\noindent Keywords: Syntactic equivalence relation; monoids, semi-rings, left skew braces; internal categories and groupoids; Mal'tsev, Gumm and protomodular categories; connectors and centralizers of equivalence relations.\\
AMSclass: 08A30, 08B05, 08B10, 18D40, 18E13. 

\vspace{3mm}\noindent Univ. Littoral C\^ote d'Opale, UR 2597, LMPA,\\
Laboratoire de Math\'ematiques Pures et Appliqu\'ees Joseph Liouville,\\
F-62100 Calais, France.\\
bourn@univ-littoral.fr

\end{document}